\documentclass[12pt]{amsart}
\usepackage{amsmath}
\usepackage{latexsym}
\usepackage{amssymb}
\usepackage{amscd}

\usepackage{amsfonts}
\usepackage[all]{xy}

\newtheorem{proposition}{Proposition}

\newtheorem{lemma}[proposition]{Lemma}
\newtheorem{main-lemma}[proposition]{Main Lemma}
\newtheorem{definition}[proposition]{Definition}
\newtheorem{lemma-definition}[proposition]{Lemma-Definition}
\newtheorem{theorem}[proposition]{Theorem}

\newtheorem{conjecture}[proposition]{Conjecture}
\newtheorem{question}[proposition]{Question}

\newtheorem{remark}[proposition]{Remark}

\newenvironment{prf}{\par\vspace{-5pt}%
\par\noindent\begingroup%
\leftskip=0em\hspace{0em}{\bf Proof.}}%
{\endgroup\hfill$\Box$}

\newcounter{tmp}

\def\db#1{ \mathbf{D}^b({#1})}
\def\perf#1{{\mathfrak P}{\mathfrak e}{\mathfrak r}{\mathfrak f}({#1})}

\def\coh{\operatorname{coh}}

\def\lto{\longrightarrow}

\def\T{{\mathcal T}}
\def\I{{\mathcal I}}
\def\O{{\mathcal O}}

\def\L{{\mathcal L}}
\def\E{{\mathcal E}}
\def\F{{\mathcal F}}
\def\G{{\mathcal G}}

\def\M{{\mathcal M}}
\def\N{{\mathcal N}}

\def\ZZ{{\mathbb Z}}

\def\NN{{\mathbb N}}
\def\ZZ{{\mathbb Z}}

\def\PP{{\mathbb P}}

\def\Hom{\operatorname{Hom}}

\def\Spec{\operatorname{Spec}}
\def\Pic{\operatorname{Pic}}

\def\Ext{\operatorname{Ext}}

\def\Ker{\operatorname{Ker}}
\def\Coker{\operatorname{Coker}}

\title[]{Remarks on generators and dimensions of triangulated categories}

\author[]{Dmitri Orlov}

\address{Algebra Section, Steklov Math. Institute of RAS, 8 Gubkin str., Moscow, 119991 Russia}

\email{orlov@mi.ras.ru}

\thanks{The author was partially supported by grant RFFI 05-01-01034, grant INTAS
05-1000008-8118, grant NSh-9969.2006.1 and Oswald Veblen Fund}

\date{}

\begin{document}
\begin{abstract}
In this paper we prove that the dimension of the bounded  derived category of coherent sheaves on a smooth
quasi-projective curve is equal to one. We also discuss dimension spectrums of these categories.
\end{abstract}

\maketitle

Let $\T$ be a triangulated category. We say that an object $E\in \T$
is a {\sf classical generator} for $\T$ if the category $\T$ coincides with the smallest triangulated subcategory of
$\T$ which contains $E$ and is closed under direct summands.

If a classical generator generates the whole category for a finite number of steps then it called
a {\sf strong generator}.
More precisely, let $\I_1$ and $\I_2$ be two full subcategories of $\T.$ We denote by $\I_1*\I_2$ the full subcategory of $\T$
consisting of all objects such that there is a distinguished triangle $M_1\to M\to M_2$ with $M_i\in \I_i.$
For any subcategory $\I\subset\T$ we denote by $\langle \I\rangle$ the smallest full subcategory of $\T$ containing $\I$ and closed under
finite direct sums, direct summands and shifts. We put $\I_1 \diamond\I_2=\langle \I_1*\I_2\rangle$ and we define by induction
$\langle \I\rangle_k=\langle\I\rangle_{k-1}\diamond\langle \I\rangle.$ If $\I$ consists of an object $E$ we denote $\langle \I\rangle$ as
 $\langle E\rangle_1$ and put by induction $\langle E\rangle _k=\langle E\rangle_{k-1}\diamond\langle E\rangle_1.$
\begin{definition}
Now we say that $E$ is a {\sf strong generator} if $\langle E\rangle_n=\T$ for some $n\in\NN.$
\end{definition}
Note that $E$ is classical generator if and only if
$\mathop\bigcup\limits_{k\in\ZZ} \langle E\rangle_k=\T.$
It is also easy to see that if a triangulated category $\T$ has a strong generator then any
classical generator of $\T$ is  strong as well.

Following to \cite{Ro} we define the dimension of a triangulated
category.
\begin{definition}
The {\sf dimension} of a triangulated category $\T,$ denoted by $\dim \T,$ is the minimal integer $d\ge 0$ such that there is $E\in \T$
with $\langle E\rangle_{d+1}=\T.$
\end{definition}
We also can define the dimension spectrum of a triangulated category as follows.
\begin{definition}
The {\sf dimension spectrum} of a triangulated category $\T,$ denoted by $\sigma(\T),$ is a subset of $\ZZ,$ which consists of
all integer $d\ge 0$ such that there is $E\in \T$
with $\langle E\rangle_{d+1}=\T$ and $\langle E\rangle_{d}\ne\T.$
\end{definition}

A.~Bondal and M.~Van den Bergh showed in \cite{BV} that the triangulated
category of perfect complexes $\perf{X}$ on a quasi-compact quasi-separated scheme $X$ has a classical generator.
(Recall that a complex of $\O_X$\!-modules is called perfect if it is locally quasi-isomorphic to a bounded complex of vector bundles.)

For the triangulated category of perfect complexes on a quasi-projective scheme we can present a classical generator directly.
\begin{theorem}
Let $X$ be a quasi-projective scheme of dimension $d$ and let $\L$ be a very ample line bundle on $X.$
Then the object $\E=\bigoplus_{i=k-d}^{k}\L^{i}$ is a classical generator for the triangulated category of perfect complexes
$\perf{X}.$
\end{theorem}
\begin{prf}
The scheme $X$ is an open subscheme of a projective scheme $X'\subset \PP^N$ and $\L$ is the restriction of $\O_{\PP^N}(1)$
on $X.$
Let us take $N+1$ linear independent hyperplanes $H_i\subset \PP^N, i=0,...,N.$ In this case
the intersection $H_0\cap\cdots\cap H_{N}$ is empty. The hyperplanes $H_i$ give a section $s$ of the vector
bundle $U=O(1)^{\oplus (N+1)}$ which does not have zeros. This implies that the Koszul complex induced by $s$
$$
0\lto\Lambda^{N+1}(U^*)\lto \Lambda^{N}(U^*)\lto\cdots\lto \Lambda^2(U^*)\lto U^{*}\lto\O_{\PP^N}\lto 0
$$
is exact on $\PP^N.$ Consider the restriction  of the truncated complex on $X$
$$
\Lambda^{d+1}(U^*_{X})\lto\cdots\lto \Lambda^2(U^*_{X})\lto U^{*}_X.
$$
It has two nontrivial cohomologies, one of which is  $\O_X.$ And, moreover, since the dimension of $X$ is equal to $d$ the sheaf $\O_X$
is a direct summand of this complex.
Tensoring this complex with $\L^{k+1}$ we obtain that the triangulated subcategory which contains $\L^{i}$ for
$i=k-d,\dots, k$ also contains $\L^{k+1}.$ Thus, it contains $\L^i$ for all $i\ge k-d.$ By duality this category contains
also all $\L^i$ for all $i\le k.$ Thus we have all powers $\L^i,$ where $i\in \ZZ.$

Finally, it easy to see that $\{\L^i\}_{i\in \ZZ}$ classically generate the triangulated category
of perfect complexes $\perf{X}.$ Indeed, for any perfect complex $E$ we can construct a bounded above complex
$P^{\cdot},$ where all $P^k$ are direct sums of line bundles $\L^i,$ together with a quasi-isomorphism  $P^{\cdot}\stackrel{\sim}{\lto}E.$
Consider the brutal truncation $\sigma^{\ge -m}P^{\cdot}$ for sufficiently large $m$ and the map $\sigma^{\ge -m}P^{\cdot}\lto E.$ The cone
of this map is isomorphic to $\F[m+1],$ where $\F$ is a vector bundle. And since the $\Hom(E, \F[m+1])=0$ for sufficiently large $m$
we get that $E$ is a direct summand of $\sigma^{\ge -m}P^{\cdot}.$
\end{prf}

A.~Bondal and M.~Van den Bergh also proved that for any smooth separated scheme $X$ the triangulated category of perfect complexes
$\perf{X}$  has a strong generator (\cite{BV}, Th.3.1.4).
Furthermore, R.~Rouquier showed that for quasi-projective scheme $X$ the property to be regular is equivalent to the property
that the triangulated category of perfect complexes $\perf{X}$ has a strong generator (see \cite{Ro}, Prop 7.35).
On the other hand, there is a remarkable result of R.~Rouquier which says that under some general conditions the
bounded derived category of coherent sheaves $\db{\coh(X)}$ has a strong generator. More precisely it says

\begin{theorem}\rm{(R.~Rouquier, \cite{Ro} Th.7.39)} Let $X$ be a separated scheme of finite type. Then there are an object
$E\in \db{\coh(X)}$ and an integer $d\in \ZZ$ such that $\db{\coh(X)}\cong\langle E\rangle_{d+1}.$ In particular,
$\dim \db{\coh(X)}<\infty.$
\end{theorem}

Keeping in mind this theorem  we can ask about the dimension of the derived category of coherent sheaves on a separated scheme of finite type.
It is proved in \cite{Ro} that
\begin{itemize}
\item for a reduced separated scheme $X$ of finite type  $\dim\db{\coh(X)}\ge\dim X;$
\item for a smooth affine scheme  $\dim\db{\coh(X)}=\dim X;$
\item for a smooth quasi-projective scheme $\dim\db{\coh(X)}\le 2\dim X.$
\end{itemize}
In this paper we show that the dimension of the derived category of coherent sheaves on a smooth quasi-projective curve $C$
is equal to 1. For affine curve it is known and for $\PP^1$ it is evident. Thus, it is sufficient
to consider a smooth projective curve of genus $g\ge 1.$

\begin{theorem}\label{th}
Let $C$ be a smooth projective curve of genus $g\ge 1.$ Then $\dim\db{\coh(C)}=1.$
\end{theorem}

At first, we should bring an object which generates $\db{\coh(C)}$ for one step. Let $\L$ be a line bundle on $C$ such that
$\deg\L\ge 8g.$ Let us consider $\E=\L^{-1}\oplus\O_C\oplus\L\oplus\L^2.$
We are going to show that $\E$ generates the bounded derived category of coherent sheaves on $C$ for one step, i.e. $\langle\E\rangle_2=\db{\coh X}.$

Since any object of $\db{\coh(X)}$ is a direct sum of its cohomologies it is sufficient to prove that
any coherent sheaf $\G$ belongs to $\langle \E\rangle_2.$ Further, each coherent sheaf $\G$ on a curve is a direct sum
of a torsion sheaf $T$ and a vector bundle $\F.$

\begin{lemma}\label{tors}Let $C$ be a smooth projective curve of genus $g\ge 1$ and
let $\L$ be a line bundle on $C$ as above. Then there is an exact sequence of the form
$$
(\L^{-1})^{\oplus r_1}\lto\O_C^{\oplus r_0}\lto T\lto 0
$$
for any torsion coherent sheaf $T$ on $C.$
\end{lemma}

Let $\F$ be a vector bundle on the curve $C.$ Consider the Harder-Narasimhan filtration
$0=\F_0\subset\F_1\subset\cdots\subset\F_n=\F.$ It is such filtration that every quotient $\F_i/\F_{i-1}$ is semi-stable and
$\mu(\F_i/\F_{i-1})>\mu(\F_{i+1}/\F_{i})$ for all $0<i<n,$ where $\mu(\G)$ is the slope of a vector bundle $\G$ and is equal to $c_1(\G)/r(\G).$

\begin{main-lemma}\label{main}
Let $\L$ be a line bundle with $\deg\L\ge 8g.$  Let $\F$ be a vector bundle on $C$ and let $0=\F_0\subset\F_1\subset\cdots\subset\F_n=\F$
be its Harder-Narasimhan filtration.
Choose $0\le i\le n$ such that $\mu(\F_i/\F_{i-1})\ge 4g>\mu(\F_{i+1}/\F_{i}).$
Then there are exact sequences of the form
$$
a)\quad (\L^{-1})^{\oplus r_1}\stackrel{\alpha}{\lto}\O_C^{\oplus r_0}\lto \F_i\lto 0,
\qquad b)\quad 0\lto \F/\F_i\lto\L^{\oplus s_0}\stackrel{\beta}{\lto}(\L^2)^{\oplus s_1}.
$$
\end{main-lemma}

To prove Lemma \ref{tors} and the Main Lemma \ref{main} we need the following lemma which is well-known.
\begin{lemma}\label{stab} Let $\G$ be a vector bundle on  a smooth projective curve $C$ over a field $k.$
Denote by $\overline{\G}$ its pullback on $\overline{C}=C\otimes_{k}\overline{k}.$
Assume that for any line bundle $\M$ on $\overline{C}$ with $\deg \M=d$ we have $H^1( \overline{C}, \overline{\G}\otimes\M)=0.$
Then
\begin{itemize}
\item[{\rm i)}] $H^1( C, \G\otimes\N)=0$ for any $\N$ on $C$ with $\deg\N\ge d;$
\item[{\rm ii)}] any sheaf $\G\otimes\N$ is generated by the global sections for all $\N$ with $\deg\N>d.$
\end{itemize}
\end{lemma}
\begin{prf}
i) Since any field extension is strictly flat it is sufficient to check that $H^1( \overline{C}, \overline{\G}\otimes\overline{\N})=0.$
From an exact sequence
\begin{equation}\label{seq}
0\to
\overline{\G}\otimes\overline{\N}(-x)\to\overline{\G}\otimes\overline{\N}\to
(\overline{\G}\otimes\overline{\N})_x\to 0
\end{equation}
on $\overline{C}$ we deduce that if $H^1(\overline{C},
\overline{\G}\otimes\overline{\N}(-x))=0$ then $H^1(\overline{C}, \overline{\G}\otimes\overline{\N})=0.$
This implies i).

ii) By the same reason as above it is enough to show that the sheaf $\overline{\G}\otimes \overline{\N}$ is generated by the global sections.
Since by  $H^1(\overline{C}, \overline{\G}\otimes \overline{\N}(-x))=0$ the map
$$
H^0(\overline{C}, \overline{\G}\otimes \overline{\N})\to
H^0(\overline{C}, (\overline{\G}\otimes \overline{\N})_x)
$$
is surjective
for any $x\in \overline{C}.$ Hence, $\overline{\G}\otimes \overline{\N}$ and $\G\otimes\N$ are generated by the global sections for all
$\N$ of degree greater than $d.$
\end{prf}

{\bf\noindent Proof of Lemma \ref{tors}.\;} Any torsion sheaf $T$ is generated by the global sections.
Consider the surjective map $\O_C^{\oplus r_0}\to T,$ where $r_0=\dim H^0(T).$ Denote by $U$ the kernel of this map.
Now it is evident that $H^1(\overline{U}\otimes\M)=0$ for any line bundle $\M$ on $\overline{C}$ with
$\deg\M\ge 2g-1,$ because $H^1(\M)=0.$ Applying Lemma \ref{stab} we get that $U\otimes\L$ is generated by the global sections.
Hence, there is an exact sequence of the form
$$
(\L^{-1})^{\oplus r_1}\lto\O_C^{\oplus r_0}\lto T\lto 0.
$$
for any torsion sheaf $T.$\hfill$\Box$

{\bf\noindent Proof of the Main Lemma.\;}
If $\G$ is a semi-stable vector bundle on $C$ with $\mu(\G)\ge 2g$ then by Serre duality we have $H^1(\overline{C}, \overline{\G}\otimes\M)=0$ for
all $\M$ with $\deg\M\ge -1.$ Therefore, by Lemma \ref{stab} the bundle $\G$ is generated by the global sections.

Now  $\F_i\subseteq\F$ as an extension of semi-stable sheaves with $\mu\ge 4g$ is generated by the global sections as well.
Consider the short exact sequence
$$
0\lto U\lto \O_C^{\oplus r_0}\lto \F_i\lto 0,
$$
where $r_0$ is the dimension of $H^0(\F_i).$ Take a line bundle $\M$
on $\overline{C}$ of degree $2g$ and consider the diagram
$$
\begin{CD}
&&0&& 0 && 0\\
&& @VVV @VVV @VVV \\
0 @>>> \overline{U}\otimes\M^{-1} @>>> {(\M^{-1})}^{\oplus r_0} @>>>  \overline{\F}_i\otimes\M^{-1} @>>> 0\\
&& @VVV @VVV @VVV \\
0 @>>> \overline{U}^{\oplus 2}@>>> \O_{\overline{C}}^{\oplus 2r_0} @>>>  \overline{\F}_i^{\oplus 2} @>>> 0\\
&& @VVV @VVV @VVV \\
0 @>>> \overline{U}\otimes\M @>>> {\M}^{\oplus r_0} @>>> \overline{\F}_i\otimes\M @>>> 0\\
&& @VVV @VVV @VVV \\
&&0&& 0 && 0
\end{CD}
$$
Since the sheaf $\overline{\F}_i\otimes\M^{-1}$ is the extension of
semi-stable sheaves with $\mu\ge 2g$ we have
$H^1(\overline{\F}_i\otimes\M^{-1})=0.$ Hence, the map
$H^0(\overline{\F}_i^{\oplus 2})\to H^0(\overline{\F}_i\otimes\M)$
is surjective. Further, we know that the map
$H^0(\O_{\overline{C}}^{\oplus 2r_0})\to H^0(\overline{\F}_i^{\oplus
2})$ is surjective and the map $H^0(\O_{\overline{C}}^{\oplus
2r_0})\to H^0(\M^{\oplus r_0})$ is injective. This implies that the
map $H^0(\M^{\oplus r_0})\to H^0(\overline{\F}_i\otimes\M)$ is
surjective as well. Hence $H^{1}(\overline{U}\otimes\M)=0.$
Therefore, by Lemma \ref{stab} the bundle $\overline{U}\otimes\M'$
is generated by the global sections for all $\M'$ with $\deg\M'\ge
2g+1.$ In particular, $U\otimes\L$ is generated by the global
sections. Thus, we get an exact sequence
$$
\quad (\L^{-1})^{\oplus r_1}\stackrel{\alpha}{\lto}\O_C^{\oplus r_0}\lto \F_i\lto 0.
$$
Sequence b) can be obtained by dualizing of  sequence a) applied for the sheaf $\F^*\otimes\L.$
\hfill$\Box$

{\bf\noindent Proof of Theorem \ref{th}.\;}
At first, since the category of coherent sheaves on $C$ has homological dimension one we see that any torsion sheaf
$T$ is a direct summand of the complex of the form $(\L^{-1})^{\oplus r_1}{\to}\O_C^{\oplus r_0}.$ Hence, it belongs to
$\langle \E\rangle_2.$

Now consider a vector bundle $\F$ on $C$ with the Harder-Narasimhan filtration $0=\F_0\subset\F_1\subset\cdots\subset\F_n=\F.$
As above let us fix $0\le i\le n$ such that $\mu(\F_i/\F_{i-1})\ge 4g>\mu(\F_{i+1}/\F_{i}).$
Applying the Main Lemma we obtain the following long exact sequence
$$
0\lto \Ker\alpha\lto (\L^{-1})^{\oplus r_1}\stackrel{\alpha}{\lto}\O_C^{\oplus r_0}\lto \F \lto\L^{\oplus s_0}\stackrel{\beta}{\lto}(\L^2)^{\oplus s_1}
\lto \Coker \beta\lto 0.
$$
Furthermore, it is easy to see that the canonical map $\Ext^1(\L^{\oplus s_0}, \O_C^{\oplus r_0})\lto
\Ext^1(\F/\F_i, \F_i)$ is surjective. Let us fix $e\in \Ext^1(\F/\F_i, \F_i)$ which defines $\F$ as the extension and choose  some
its pull back  $e'\in \Ext^1(\L^{\oplus s_0}, \O_C^{\oplus r_0}).$

Now let us consider the map
\begin{equation}\label{cone}
\phi: (\L^{-1})^{\oplus r_0}\oplus\L^{\oplus s_0}[-1]\lto \O_C^{\oplus r_0}\oplus(\L^2)^{\oplus s_1}[-1],
\quad\text{
where}\quad \phi=
\begin{pmatrix}
\alpha & e'\\
0 & \beta
\end{pmatrix}
\end{equation}
and take a cone $C(\phi)$ of $\phi.$ The cone $C(\phi)$ is isomorphic to a complex that has  three nontrivial cohomologies
$H^{-1}(C(\phi))\cong \Ker\alpha,\; H^{1}(C(\phi))\cong\Coker\beta$ and, finally, $H^0(C(\phi))\cong\F.$
Thus, $\F$ is a direct summand of $C(\phi)$ and, consequently, it belongs to $\langle\E\rangle_2.$ This implies
that the whole bounded derived category of coherent sheaves on $C$ coincides with $\langle \E\rangle_2$ and the dimension of
$\db{\coh(C)}$ is equal to 1.
\hfill$\Box$

Having in view of the given theorem we may assume, that the following conjecture can be true.
\begin{conjecture}
Let $X$ be a smooth quasi-projective scheme of dimension $n.$ Then $\dim\db{\coh(X)}=n.$
\end{conjecture}
\begin{remark}{\rm
For a non regular scheme it is evidently not true. For example, the dimension of the bounded derived
category of coherent sheaves on the zero-dimension scheme $\Spec(k[x]/x^2)$ equals to 1.}
\end{remark}

It is also very interesting to understand what the  spectrum $\sigma(\db{\coh(X)})$ forms.
In particular we can ask the following questions
\begin{question}
Is the spectrum of the bounded derived category of coherent sheaves on a smooth quasi-projective scheme
bounded? Is it bounded for a non smooth scheme?
\end{question}
\begin{question}
Does the spectrum of the bounded derived category of coherent sheaves on a (smooth) quasi-projective scheme
form an integer interval?
\end{question}
Let us try to calculate  the dimension spectra of the derived categories of coherent sheaves on some smooth curves.
\begin{proposition}
Let $C$ be a smooth affine curve. Then the dimension spectrum $\sigma(\db{\coh C})$ coincides with
$\{1\}.$
\end{proposition}
\begin{prf}
If $\E$ is a strong generator then it has a some locally free sheaf $\F$ as a direct summand.
Now since $C$ is affine then there is an exact sequence of the form
$$
\F^{\oplus r_1}\lto \F^{\oplus r_0}\lto \G\lto 0
$$
for any coherent sheaf $\G$ on $C.$ Hence, any coherent sheaf $\G$ belongs to $\langle\E\rangle_2.$
Since the global dimension of $\coh C$ is equal to 1 we obtain that $\langle\E\rangle_2=\db{\coh C}.$
\end{prf}

We can also find the dimension spectrum of the projective line.
\begin{proposition}
The dimension spectrum $\sigma(\db{\coh\PP^1})$ coincides with the set $\{1,2\}.$
\end{proposition}
\begin{prf}
Indeed, 1 is the dimension. And, for example, the object $\E=\O(-1)\oplus\O$ generate the whole category
$\db{\coh\PP^1}$ for one step. Now, the object $\E=\O_{\PP^1}\oplus \O_p,$ where $p$ is a point, is a generator, because
$\O(-1)$ belongs to $\langle \E\rangle_2.$ This also implies that $\langle\E\rangle_3\cong\db{\coh\PP^1}.$
On the other hand, $\langle \E\rangle_2\ncong\db{\coh\PP^1}.$ To see it we can check that an object $\O_q,$ where $q\ne p,$ doesn't belong
to $\langle \E\rangle_2.$ Indeed, $\O_q$ is completely orthogonal to $\O_p$ and doesn't belong to subcategory generated by $\O.$
Finally, it easy to see that any object $\E,$ which  generates the whole category, generates it at least for two steps, i.e.
$\langle \E\rangle_3\cong\db{\coh\PP^1}.$ If $\E$ contains as direct summands two different line bundles than it generates
the whole category for one step. If $\E$ has only one line bundle as a direct summand then it also has a torsion sheaf as a direct summand.
This implies that $\langle\E\rangle_2$ has another line bundle. Therefore, $\langle \E\rangle_3$ is the whole category.
\end{prf}

Another simple result says
\begin{proposition}
Let $C$ be a smooth projective curve of genus $g>0$ over a field $k.$ Assume that $C$ has at least two different points
over $k.$ Then the dimension spectrum $\sigma(\db{\coh C})$ contains $\{1,2\}$ as a proper subset,
i.e. $\{1,2\}$ is strictly contained in the dimension spectrum.
\end{proposition}
\begin{prf}
The spectrum contains 1 as the dimension of the category. Let us now take a line bundle $\L$ on $C$ which satisfies the condition
as in Theorem \ref{main}, i.e. $\deg \L\ge 8g$ and consider the object $\E=\O_C\oplus \L^2.$
It is easy to see that the line bundles $\L^{-1}$ and $\L$ belong to $\langle \E\rangle_2,$  because there are exact sequences
$$
0\lto \O_C\lto \L^{\oplus 2}\lto \L^2\lto 0\qquad\text{and}\qquad
0\lto (\L^{-1})^{\oplus 2}\lto \O_C^{\oplus 3}\lto \L^{2}\lto 0.
$$
The proof of Theorem \ref{main} (see the map \ref{cone}) implies that $\langle\E\rangle_3\cong\db{\coh C}.$
On the other hand, the subcategory $\langle\E\rangle_2$ doesn't coincides with the whole $\db{\coh C}.$
For example, a nontrivial line bundle $\M$ from $\Pic^{0} C$ doesn't belong to $\langle \E\rangle_2,$ because it is
completely orthogonal to the structure sheaf $\O_C$ and, evidently, could not be obtained from the line bundle $\L^2.$

Let us take a point $p\in C$ and consider the object $\E=\O_C\oplus\O_p,$ where $\O_p$ is the skyscraper in $p.$
This object is a strong generator and we can show that $\langle\E\rangle_3\ne \db{\coh C}.$ Take another point $q\ne p$
and consider the skyscraper sheaf $\O_q.$ It is completely orthogonal to $\O_p$ and have only one-dimensional 1-st Ext
to $\O_C.$ Hence, if $\O_q$ belongs to $\langle\E\rangle_3$ then it should be a direct summand of an object $M$ which
is included in an exact triangle of the form
$$
\O_C^{\oplus k}\lto N\lto M\lto \O_C^{\oplus k}[1],
$$
where $N\in\langle\E\rangle_2.$ Since the 1-st Ext from $\O_q$ to $\O_C$ is one-dimensional we can take $k=1.$ The composition
of the map $\O_q\to M$ with $M\to \O_C$ should be the nontrivial 1-st Ext from $\O_q$ to $\O_C.$
Now object $N$ is a direct sum  of indecomposable objects from $\langle\E\rangle_2.$ It is easy to see that we can consider
only objects for which there are nontrivial homomorphisms from $\O_C$ and nontrivial homomorphisms to $\O_q.$ All other can be removed from $N.$
Thus $N$ is a direct sum of $\O(p)$ and objects $U$ that are extensions
\begin{equation}\label{ext}
0\lto \O_C^{\oplus r_1}\lto U\lto \O_C^{\oplus r_2}\lto 0.
\end{equation}
Finally, split embedding $\O_q\to M$ gives us a nontrivial map from $\O(q)$ to $N.$ But there are no nontrivial maps from $\O(q)$ to
$\O(p)$ and to $U$ of the form (\ref{ext}). Therefore, $\O_q$ can not belong to $\langle\E\rangle_3.$
\end{prf}

\noindent{\bf Acknowledgments:} I am grateful to Igor Burban and
Pierre Deligne for useful notes and to Calin Lazaroiu, Tony Pantev
and Raphael Rouquier for very interesting discussions. I would like
to thank Institute for Advanced Study, Princeton, USA for a
hospitality and a very stimulating atmosphere.


\end{document}